\documentclass[12pt,reqno]{amsart}
\usepackage{graphics}
\usepackage{graphicx} 
\usepackage{psfrag} 
\usepackage{amssymb,a4wide,amsthm,amsmath,amsfonts}
\usepackage{tikz}
\usetikzlibrary{matrix}
\usepackage{bbm}
\usepackage{dsfont}
\usepackage{mathrsfs}
\usepackage{color}

\theoremstyle{plain}
\newtheorem{theorem}{Theorem}[section]
\theoremstyle{remark}

\theoremstyle{plain}
\newtheorem{corollary}[theorem]{Corollary}
\newtheorem{lemma}[theorem]{Lemma}
\newtheorem{proposition}[theorem]{Proposition}

\numberwithin{equation}{section}

\newcommand{\R}{\mathbb{R}}
\newcommand{\N}{\mathbb{N}}
\newcommand{\Z}{\mathbb{Z}}

\newcommand{\F}{\mathcal{F}}

\newcounter{probcount}

\begin{document}
\title{Products of Lipschitz-free spaces and applications}
\author[P. L. Kaufmann]{Pedro L. Kaufmann}
\address{CAPES Foundation, Ministry of Education of Brazil, Bras\'ilia/DF 70040-020, Brazil and Institute de Math\'ematiques de Jussieu\\ Universit\'e Pierre et Marie Curie\\ 4 Place Jussieu, 75005 Paris, France}

\begin{abstract}

We show that, given a Banach space $X$, the Lipschitz-free space over $X$, denoted by $\mathcal{F}(X)$, is isomorphic to $(\sum_{n=1}^\infty \mathcal{F}(X))_{\ell_1}$. Some applications are presented, including a non-linear version of Pe\l czy\'ski's decomposition method for Lipschitz-free spaces and  the identification up to isomorphism between $\mathcal{F}(\mathbb{R}^n)$ and the Lipschitz-free space over any compact metric space which is locally bi-Lipschitz embeddable into $\mathbb{R}^n$ and which contains a subset that is Lipschitz equivalent to the unit ball of $\mathbb{R}^n$. We also show that $\mathcal{F}(M)$ is isomorphic to $\mathcal{F}(c_0)$ for all separable metric spaces $M$ which are absolute Lipschitz retracts and contain a subset which is Lipschitz equivalent to the unit ball of  $c_0$. This class contains all $C(K)$ spaces with $K$ infinite compact metric (Dutrieux and Ferenczi had already proved that $\mathcal{F}(C(K))$ is isomorphic to $\mathcal{F}(c_0)$ for those $K$ using a different method). 
Finally we study Lipschitz-free spaces over  certain unions and quotients of metric spaces, extending a result by Godard.

\end{abstract}



\maketitle

\section{Introduction} 


Let $(M,d,0)$ be a pointed metric space (that is, a distinguished point $0$ in $M$, called \emph{base point}, is chosen), and consider the Banach space $Lip_0(M)$ of all real-valued Lipschitz functions on $M$ which vanish at $0$, equipped with the norm 
$$\|f\|_{Lip}:= \inf_{x,y\in M,x\neq y} \frac{|f(x)-f(y)|}{d(x,y)}.$$
On the closed unit ball of $Lip_0(M)$, the topology of pointwise convergence is compact, so $Lip_0(M)$ admits a canonical predual, which is called the \emph{Lipschitz-free space over $M$} and denoted by $\F(M)$. This space is the closure in $Lip_0(M)^*$ of $span\{\delta_x:x\in M\}$, where $\delta_x$ are the evaluation functionals defined by $\delta_x(f)=f(x)$. It is readily verified that $\delta:x\mapsto \delta_x$ is an isometry from $M$ into $\F(M)$.  Given $0'$, it is clear that $T:Lip_0(M)\rightarrow Lip_{0'}(M)$ defined by $T(f):=f-f(0')$  is a weak$^*$-to-weak$^*$ continuous isometric isomorphism, thus the choice of different base points yields isometrically isomorphic Lipschitz-free spaces. We refer to \cite{W} for a study  on Lipschitz functions spaces, and to \cite{W} and \cite{GK} for an introduction to Lipschitz-free spaces and its basic properties. 

One of the main properties of the Lipschitz-free spaces is that it permits to interpret Lipschitz maps between metric spaces from a linear point of view:

\begin{proposition}

Let $M$ and $N$ be pointed metric spaces, consider $\delta^M$ and $\delta^N$ the isometries that assign each $x\in M$ (respectively, each $x\in N$) to the corresponding evaluation functional $\delta^M_x$ in $\F(M)$ (respectively, $\delta^N_x$ in $\F(N)$), and suppose that $L:M\rightarrow N$ is a Lipschitz function such that $L(0_M)=0_N$. Then there is an unique linear map $\hat{L}:\F(M)\rightarrow \F(N)$ such that $\hat{L}\circ \delta^M = \delta^N \circ L$, that is, such that the following diagram commutes:  

\center{
\begin{tikzpicture}
  \matrix (m) [matrix of math nodes,row sep=3em,column sep=4em,minimum width=2em]
  {
     M & N \\
     \F(M) & \F(N) \\};
  \path[-stealth]
    (m-1-1) edge node [left] {$\delta^M$} (m-2-1)
            edge node [above] {$L$} (m-1-2)
    (m-2-1.east|-m-2-2) edge [dashed] node [above] {$\hat{L}$}
             (m-2-2)
    (m-1-2) edge node [right] {$\delta^N$} (m-2-2);
\end{tikzpicture}}

\leftline{Moreover, $\|\hat{L}\|=\|L\|_{Lip}$.}
\label{propgk}
\end{proposition}

In particular, if $M$ and $N$ are Lipschitz equivalent (that is, there is a bi-Lipschitz map between $M$ and $N$) then $\F(M)$ and $\F(N)$ are isomorphic. The converse is not true, even if $M$ and $N$ are assumed to be Banach spaces: if $K$ is an infinite compact metric space, then $\F(C(K))$ is isomorphic to $\F(c_0)$, even though $C(K)$ is not Lipschitz equivalent to $c_0$ in general (recall that, if $C(K)$ is uniformly homeomorphic to $c_0$, then it is isomorphic to $c_0$ c.f. \cite{JLS}). This first counterexample for the Banach space case was presented by Dutrieux and Ferenczi in \cite{DF}. 

Despite of the simplicity of the definition of the Lipschitz-free spaces, many fundamental questions about their structure remain unanswered. Godard \cite{G} characterized the metric spaces $M$ such that $\F(M)$ is isometrically isomorphic to a subspace of $L^1$ as exactly those who are isometrically embeddable into $\R$-trees (that is, connected graphs with no cycles with the graph distance); on the other hand, Naor and Schechtman \cite{NS} have shown that $\F(\Z^2)$ (thus also  $\F(\R^2)$) is not isomorphic to any subspace of $L^1$. From this arises the natural question of characterizing the metric spaces $M$ such that $\F(M)$ is (nonisometrically) isomorphic to $L^1$.  Godefroy and Kalton \cite{GK} have shown that, given a Banach space $X$, $X$ has the bounded approximation property if and only if $\F(X)$ has that property, and recently H\'ajek and Perneck\'a \cite{HP} have shown that $\F(\R^n)$ admits a Schauder basis, refining a result from \cite{LP}. However, it is not known whether $\F(F)$ admits a Schauder basis for any given closed subset $F\subset \R^n$. It is also not of this author's knowledge whether $\F(\R^n)$ is isomorphic to $\F(\R^m)$ or not, for distinct $m,n\geq 2$. Even the study of Lipschitz-free spaces over very simple subsets of $\R^2$ can present difficulties (see the question posed after Proposition \ref{propbasic2}). 

In this context, we continue the exploration of what could be considered  basic properties of Lipschitz-free spaces and their relation with the underlying metric spaces. We will show, for instance, that for any given Banach space $X$, we have that $\F(X)$ is isomorphic to $(\sum_{n=1}^\infty \F(X))_{\ell_1}$ (see Theorem \ref{main1}). This provides in particular a kind of non-linear version for Pe\l czy\'ski's decomposition method (see Corollary \ref{corolpeln}), which in turn can be used to obtain the mentioned example by Dutrieux and Ferenczi of non-Lipschitz equivalent Banach spaces sharing the same Lipschitz-free space. In fact, we show that $\F(M)$ is isomorphic to $\F(c_0)$ for a wider class of metric spaces  (see Corollary \ref{corolfxfc0}). We will also show that, for compact metric spaces $M$ which are locally bi-Lipschitz embeddable in $\R^n$, we have that $\F(M)$ admits a complemented copy in $\F(\R^n)$; when moreover the euclidean ball $B_{\R^n}$ is bi-Lipschitz embeddable in $M$, we have that $\F(M)$ and $\F(\R^n)$ are actually isomorphic (see Theorem \ref{main2}). The class of metric spaces satisfying both properties includes all $n$-dimensional compact Riemannian manifolds. 
Independently, we study the behavior of Lipschitz-free spaces with respect to certain gluings of metric spaces, expanding an initial idea presented by Godard \cite{G}. 

\subsection{Notation} We say that two metric spaces $M$ and $N$ are \emph{$C$-Lipschitz equivalent}, for some constant $C>0$, if there a bi-Lipschitz onto map $\varphi:M\rightarrow N$ satisfying $\|\varphi\|_{Lip}.\|\varphi^{-1}\|_{Lip}\leq C$. $M$ and $N$ are then Lipschitz equivalent if they are $C$-Lipschitz equivalent for some $C>0$; in that case we also write $M\stackrel{L}\sim N$. Given two Banach spaces $X$ and $Y$, we write $X\cong Y$ when $X$ and $Y$ are isometrically isomorphic, $X\stackrel{c}{\hookrightarrow} Y$ when there is a complemented copy of $X$ in $Y$, and $X\simeq Y$ when $X$ and $Y$ are isomorphic. If $X$ and $Y$ are isomorphic, the \emph{Banach-Mazur distance between $X$ and $Y$} is defined by
$$
d_{BM}(X,Y):=\inf\{\|T\|.\|T^{-1}\|:T\mbox{ is an isomorphism from $X$ onto $Y$}\}.$$  
$\|T\|.\|T^{-1}\|$ is called the \emph{(linear) distortion of $T$}. When $d_{BM}(X,Y)\leq C$ for some $C>0$, we say that \emph{$X$ is isomorphic to $Y$ with distortion bounded by $C$}. 

$Ext_0(F,M)$ denotes the set of linear extension operators for Lipschitz functions and $Ext_0^{pt}(F,M)$ is the set of pointwise-to-pointwise continuous elements of $Ext_0(F,M)$ (see Subection \ref{subseclipext}).

\subsection{Structure of this work} In Section \ref{secext}, we present some background results on linear extension operators for Lipschitz functions and consequent ways to decompose the Lipschitz-free space over a metric space using metric quotients. In Section \ref{secprod}  we show that, for every Banach space $X$, $\F(X)\simeq (\sum_{n=1}^\infty \F(X))_{\ell_1}$, and derive some consequences. 

Sections \ref{secfx2} and \ref{secunion} are independent of the results presented in Section \ref{secprod}. In Section \ref{secfx2} we show that, for every Banach space $X$, $d_{BM}(\F(X),\F(X)\oplus_1 \F(X))\leq 4$. In Section \ref{secunion} we provide formulas for Lipschitz-free spaces over certain unions of metric spaces.  \\


\section{Linear extensions of Lipschitz functions and the Lipschitz-free space over metric quotients} 
\label{secext}

\subsection{Linear extensions of Lipschitz functions}
\label{subseclipext}

Given a pointed metric space $(M,d,0)$ and a subset $F$ containing $0$, let us denote by $ Ext_0(F,M)$ the set of all extensions $E:Lip_0(F)\rightarrow Lip_0(M)$ which are linear and continuous ($E$ being an extension means that $E(f)|_F=f$ for all $f\in Lip_0(F)$). It is immediate to see that, if we choose another base point $0'$ contained in $F$, to each element $E\in Ext_0(F,M)$ there is a corresponding $E'\in Ext_{0'}(F,M)$, defined by $E'(f):=E(f-f(0'))+f(0')$, which satisfies $\|E'\|=\|E\|$, so generally it is not important which base point is chosen.  Recall that there are always continuous but not necessarily linear extensions from $Lip_0(F)$ to $Lip_0(M)$; for example the infimum convolution 
$$
E(f)(x):=\inf_{y\in F} \{f(y) + \|f\|_{Lip}d(x,y)\}
$$
is such an extension and it is an isometry, although in most cases it fails to be linear. It is possible, though, to have $ Ext_0(F,M)=\emptyset$; Brudnyi and Brudnyi provide us with an example of a two-dimensional Riemannian manifold $M$, equipped with its geodesic metric, which admits a subset $F$ satisfying that condition (see Theorem 2.18 in \cite{BB}).

We will be particularly interested in the subset  $Ext_0^{pt}(F,M)$ of $ Ext_0(F,M)$ consisting of the pointwise-to-pointwise continuous elements. The fact that on bounded sets of $Lip_0(F)$ the weak$^*$ and the pointwise topologies coincide implies that any element of $ Ext_0(F,M)$ is weak$^*$-to-weak$^*$ continuous if and only if it belongs to $Ext_0^{pt}(F,M)$. Therefore, any $E\in  Ext_0^{pt}(F,M)$ admits a preadjoint $P:\F(M)\rightarrow\F(F)$, which is a (continuous) canonical projection, in the sense that $P(\mu)=\mu|_F$ for all finitely supported $\mu\in \F(M)$. In particular, $\F(F)$ is complemented in $\F(M)$. Reciprocally, given a continuous projection $P:\F(M)\rightarrow\F(F)$ such that $P(\mu)=\mu|_F$ for all finitely supported $\mu\in \F(M)$, we have that $P^*\in  Ext_0^{pt}(F,M)$. 

Even in the context where $M$ is a Banach space and $F$ is a closed linear subspace, we might not get this complementability condition. Consider, for example, $c_0$ and let $X$ be a subspace of $c_0$ which fails to have the bounded approximation property. As we mentioned in the introduction, given a Banach space $Y$, $Y$ has the bounded approximation property if and only if $\F(Y)$ has that property. Since this property is inherited by complemented subspaces, it follows that $\F(X)$ cannot be isomorphic to a complemented subspace  of $\F(c_0)$. One can still pose the question of wether or not  $ Ext_0(X,c_0)$ is empty.

On the other hand, we have the following positive example: 

\begin{proposition}[Lancien, Perneck\'a \cite{LP}] 

There exists $C>0$ such that, for each $n\in\N$ and each subset $F$ of $\R^n$ containing $0$, there exists $E$ in $ Ext_0^{pt}(F,\R^n)$ satisfying $\|E\|\leq C\sqrt n$. 

\label{propLP}
\end{proposition}

This result appears as part of the proof of Proposition 2.3 of \cite{LP}, which states that $\F(F)$ has the $C\sqrt n$-bounded approximation property. It involves a construction by Lee and Naor from \cite{LN}, and the fact that $\R^n$ admits a so-called \emph{$K$-gentle partition of the unity with respect to $F$}, which in turn induces the mentioned extension $E$.

\subsection{Metric quotients and Lipschitz-free spaces}

We turn our attention to a special kind of metric quotient. Given a pointed metric space $(M,d,0)$ and a subset $F$ of $M$ containing $0$, let $\sim$ be the equivalence relation which collapses $\overline F$ to a point (that is, the equivalence classes are either singletons or $\overline F$). We define the metric quotient of $M$ by $F$, denoted by $M/F$, as the pointed metric space $(M/\sim,\tilde{d},[0])$, where $\tilde{d}$ is defined by
\begin{eqnarray}
\tilde{d}([x],[y])= \min\left\{d(x,y),d(x,F)+d(y,F)\right\}.
\label{deftildedF}
\end{eqnarray}
The space $Lip_{[0]}(M/F)$ can be interpreted as the closed linear subspace of $Lip_0(M)$ consisting of all of its functions which are null in $F$. Depending on how $F$ is placed in $M$, we can have the following decomposition for $\F(M)$:

\begin{lemma}

Let $(M,d,0)$ be a pointed metric space and $F$ be a subset containing $0$, and suppose that there exists $E\in Ext_0^{pt}(F,M)$. Then 
$$\F(M)\simeq \F(F)\oplus_1 \F(M/F),$$ 
with distortion bounded by $(\|E\|+1)^2$.
 
\label{lemmaextFM}
\end{lemma}

\textbf{Proof.} Define $\Phi: Lip_0(F) \oplus_\infty Lip_0(M/F)\rightarrow Lip_0(M)$  by $\Phi(f,g)\doteq E(f)+g$. It is straightforward that $\Phi$ is an onto isomorphism with $\|\Phi\|\leq \|E\|+1$, that $\Phi$ is pointwise-to-pointwise continuous and that its inverse $\Phi^{-1}: h\mapsto (h|_F,h- E(h|_F))$ has norm also bounded by $\|E\|+1$. It follows that $\Phi$ is the adjoint of an isomorphism $\Psi$ between $\F(M)$ and $\F(F)\oplus_1 \F(M/F)$ satisfying the desired distortion bound.  $\Box$\\


\section{Products of Lipschitz-free spaces}
\label{secprod}

In this section we will show that $\F(X)\simeq\left(\sum_{n=1}^\infty \F(X)\right)_{\ell_1}$ for any Banach space $X$, and derive some consequences. With that purpose we will use the following construction by Kalton \cite{K2}. Let $(M,d,0)$ be a pointed metric space, denote by $B_r$ the closed balls centered at $0$ and with radius $r>0$ and consider, for each $k\in \Z$, the linear operator $T_k:\F(M)\rightarrow \F(B_{2^{k+1}} \setminus B_{2^{k-1}})$ defined by 
$$
T_k\delta_x :=\left\{ 
\begin{array}{ll}
0,& \mbox{if } x\in B_{2^{k-1}};\\
(\log_2d(x,0)-k+1)\delta_x,& \mbox{if } x\in B_{2^{k}} \setminus B_{2^{k-1}};\\
(k+1-\log_2d(x,0))\delta_x,& \mbox{if } x\in B_{2^{k+1}} \setminus B_{2^{k}};\\
0,& \mbox{if } x\not\in B_{2^{k+1}}.
\end{array}\right.
$$ 
Lemma 4.2 from \cite{K2} says that for each $\gamma\in \F(M)$ we have that $\gamma = \sum_{k\in \Z} T_z\gamma$ unconditionally and
\begin{eqnarray}
\sum_{k\in \Z} \|T_k\gamma\|_\F \leq 72 \|\gamma\|_\F.
\label{vcvc}
\end{eqnarray}
Another result from that same paper that we will use is Lemma 4.2, which states that, given $r_1,\dots,r_n,s_1,\dots,s_n\in\Z,r_1<s_1<r_2<\dots <s_n$ and $\gamma_k\in \F(B_{2^{s_k}}\setminus B_{2^{r_k}})$ and writing $\theta :=\min_{k=1,...,n-1}\{r_{k+1}-s_k\}$, then 
\begin{eqnarray}
\|\gamma_1+\dots+\gamma_n\|_\F\geq \frac{2^\theta-1}{2^\theta+1}\sum_{k=1}^n \|\gamma_k\|_\F.
\label{wert}
\end{eqnarray}

\begin{theorem}

Let $X$ be a Banach space. Then 
$$\F(X)\simeq\left(\sum_{n=1}^\infty \F(X)\right)_{\ell_1}.$$ 

\label{main1}
\end{theorem}

\textbf{Proof.} Note that $S:(\gamma_k)\in (\sum_{k\in\Z}\F(B_{2^{k+1}} \setminus B_{2^{k-1}}))_{\ell_1}\mapsto \sum_{k\in\Z} \gamma_k \in \F(X)$ is linear, continuous and onto, and from (\ref{vcvc}) we get that $T:\gamma\in\F(X) \mapsto (T_k\gamma)\in (\sum_{k\in\Z}\F(B_{2^{k+1}} \setminus B_{2^{k-1}}))_{\ell_1}$ is a well defined one-to-one linear continuous operator. Thus $T\circ S$ is a continuous projection from $(\sum_{k\in\Z}\F(B_{2^{k+1}} \setminus B_{2^{k-1}}))_{\ell_1}$ onto the isomorphic copy $T(\F(X))$ of $\F(X)$. 

Denote $M:=\cup_{k\in\Z} (B_{2^{2k+1}}\setminus B_{2^{2k}})$, and consider $E\in  Ext_0(M\cup \{0\},X)$ which extends each element of $Lip_0(M\cup\{0\})$ linearly on each radial segment $[2^{2k-1},2^{2k}]x,k\in\Z,x\in S_X$. Clearly $E$ is pointwise-to-pointwise continuous, thus it is the adjoint of some  $P:\F(X)\rightarrow \F(M\cup\{0\})$, which is a projection which satisfies $P(\mu)=\mu|_M$ for all finitely supported $\mu\in\F(X)$. Note that $\F(M\cup\{0\})\cong \F(M)$, since $0\in\overline{M}$. Now by (\ref{wert}), the natural identification $Id:\F(M)\rightarrow (\sum_{k\in\Z} \F(B_{2^{2k+1}}\setminus B_{2^{2k}}))_{\ell_1}$ is an isomorphism. So there is a complemented copy of $(\sum_{k\in\Z} \F(B_{2^{2k+1}}\setminus B_{2^{2k}}))_{\ell_1}$ in $\F(X)$.

Note that, by Proposition \ref{propgk}, re-scalings of any metric space give rise to isometrically isomorphic Lipschitz-free spaces. Thus all spaces $\F(B_{2^{2k+1}}\setminus B_{2^{2k}}),\,k\in\Z$ are isometrically isomorphic to 
$\F(B_{2}\setminus B_{1})$ and all spaces $\F(B_{2^{k+1}} \setminus B_{2^{k-1}}),\,k\in\Z$ are isometrically isomorphic to $\F(B_4 \setminus B_1)$, which in turn is isomorphic to $\F(B_2\setminus B_1)$. It follows that 
$$
\F(X)\stackrel{c}{\hookrightarrow} \left(\sum_{j=1}^\infty \F(B_2\setminus B_1)\right)_{\ell_1}\mbox{ and } \left(\sum_{j=1}^\infty \F(B_2\setminus B_1)\right)_{\ell_1} \stackrel{c}{\hookrightarrow} \F(X). 
$$
Since $\left(\sum_{j=1}^\infty \F(B_2\setminus B_1)\right)_{\ell_1}$ is isomorphic to its $\ell_1$-sum, by a standard Pe\l czy\'nski's decomposition method we have 
\begin{eqnarray}
\F(X)\simeq \left(\sum_{j=1}^\infty \F(B_2\setminus B_1)\right)_{\ell_1}
\label{werg}
\end{eqnarray}
and the conclusion follows immediately.  $\Box$\\

As a direct consequence of Theorem \ref{main1} and Proposition \ref{propgk}, we get the following nonlinear version of Pe\l czy\'nski's decomposition method for Lipschitz-free spaces.

\begin{corollary}

Let $X$ be a Banach space and $M$ be a metric space, and suppose that $X$ and $M$ admit Lipschitz retracts $N_1$ and $N_2$, respectively, such that $X$ is Lipschitz equivalent to $N_2$ and $M$ is Lipschitz equivalent to $N_1$. Then $\F(X)\simeq \F(M)$.

\label{corolpeln}
\end{corollary}

\textbf{Proof.} $\F(X)$ is isomorphic to $\F(N_2)$, which in turn is a complemented subspace of $\F(M)$. Analogously, $\F(M)$ is isomorphic to a complemented subspace of $\F(N_1)$. The conclusion follows by applying the standard Pe\l czy\'nski's decomposition method.  $\Box$\\

\begin{corollary}

Let $X$ be a Banach space. Then 
$$
\F(X)\simeq \F(B_1).
$$

\label{propfxfb}
\end{corollary}

\textbf{Proof.} Since $B_1$ is a Lipschitz retract of $X$, it follows that $\F(X)$ contains a complemented copy of $\F(B_1)$. In the proof of Theorem \ref{main1} we have shown that $\F(X)$ is isomorphic to $(\sum_{k\in\Z} \F(B_{2^{2k+1}}\setminus B_{2^{2k}}))_{\ell_1}$, which is clearly isomorphic to 
 $(\sum_{k<0} \F(B_{2^{2k+1}}\setminus B_{2^{2k}}))_{\ell_1}$ since all summands are isometrically isomorphic. Let $N:=\cup_{k<0} (B_{2^{2k+1}}\setminus B_{2^{2k}})$ Again by (\ref{wert}), $(\sum_{k<0} \F(B_{2^{2k+1}}\setminus B_{2^{2k}}))_{\ell_1}$ is isomorphic to $\F(N)$, which is complemented in $\F(B_1)$ since there is a pointwise-to-pointwise continuous element in $ Ext_0(N\cup\{0\},B_1)$. The conclusion follows by an application of Pe\l czy\'nski's decomposition method.   $\Box$\\

Recall that a subset $F$ of a metric space $M$ is called a \emph{Lipschitz retract of $M$} if there is a Lipschitz map from $M$ onto $F$ which coincides with the identity on $F$;  in case such a map exists, it is called a \emph{Lipschitz retraction}. A metric space is said to be an \emph{absolute Lipschitz retract} if it is a Lipschitz retract of every metric space containing it. Given any metric space $M$, the space $C_u(M)$ of real-valued bounded and uniformly continuous functions on $M$, equipped with the uniform norm, is an example of Banach space which is an absolute Lipschitz retract (see e.g. \cite{BL}, Theorem 1.6). This class includes all $C(K)$ spaces for $K$ compact metric space, in particular it includes $c_0$. Since all separable metric spaces are bi-Lipschitz embeddable in $c_0$ (\cite{BL}, Theorem 7.11), we obtain the following  class of metric spaces $M$ with $\F(M)\simeq \F(c_0)$: 

\begin{corollary}

Let $M$ be a separable metric space containing a Lipschitz retract which is Lipschitz equivalent to the unit ball of $c_0$, and suppose that $M$ is an absolute Lipschitz retract. Then $\F(M)\simeq \F(c_0)$. 

In particular, if $K$ is an infinite compact metric space, then $\F(C(K))\simeq \F(c_0)$. 

\label{corolfxfc0}
\end{corollary}

\textbf{Proof.}  It is straightforward by Proposition \ref{propgk} and Corollary \ref{propfxfb} that there is a complemented copy of $\F(c_0)$ in $\F(M)$. $M$ is Lipschitz equivalent to some subset $F$ of $c_0$, and $F$ is an absolute Lipschitz retract since this property is preserved by Lipschitz equivalences. Thus $F$ is a Lipschitz retract of $c_0$, and again by Proposition \ref{propgk} this implies that $\F(F)$ (and thus $\F(M)$) admits a complemented copy in $\F(c_0)$. The conclusion follows from Theorem \ref{main1} and an application of Pe\l czy\'nski's decomposition method.  $\Box$\\

\begin{corollary}

Let $F$ be a subset of $\R^n$ with nonempty interior. Then $\F(F)\simeq \F(\R^n)$. 

\label{corolnonemptyint}
\end{corollary}

\textbf{Proof. } By Proposition \ref{propLP}, there is a complemented copy of $\F(F)$ in $\F(\R^n)$. Taking any closed ball $B$, it is easy to see that there is a Lipschitz retraction from $F$ onto $B$; thus by Proposition \ref{propgk} and Corollary \ref{propfxfb} there is also a complemented copy of $\F(\R^n)$ in $\F(F)$. The result follows from Theorem \ref{main1} and an application of Pe\l czy\'nski's decomposition method. .  $\Box$\\

\textbf{Remark.} In \cite{HP}, H\'ajek and Perneck\'a have shown that $\F(\R^n)$ admits a Schauder basis, and rose the natural question of whether or not the same holds true for $\F(F)$, being $F$ any closed subset of $\R^n$. Note that, by Corollary \ref{corolnonemptyint}, the problem is reduced to the case where $F$ has empty interior. \\

In order to study Lipschitz-free spaces of locally euclidean metric spaces, alongside the corollaries of Theorem \ref{main1}, the following result becomes handy:

\begin{theorem}[Lang, Plaut \cite{LanPla}]

Let $M$ be a compact metric space such that each point of $M$ admits a neighborhood which is bi-Lipschitz embeddable in $\R^n$. Then $M$ is bi-Lipschitz embeddable in $\R^n$. 

\end{theorem}

\begin{theorem}

Let $M$ be a compact metric space such that  each $x\in M$ admits a neighborhood which is bi-Lipschitz embeddable in $\R^n$. Then there is a complemented copy of  $\F(M)$ in $\F(\R^n)$. 

If moreover the unit ball of $\R^n$ is bi-Lipschitz embeddable into $M$, then $\F(M)\simeq \F(\R^n)$. In particular, the Lipschitz-free space over any $n$-dimensional compact Riemannian manifold equipped with its geodesic metric is isomorphic to  $\F(\R^n)$.

\label{main2}
\end{theorem}

\textbf{Proof.} The first part follows directly from Lang and Plaut's result and the fact that the Lipschitz-free space over any subset of $\R^n$ admits a complemented copy in $\F(\R^n)$. 

For the second part, note that the closed unit ball of $\R^n$ is an absolute Lipschitz retract, and recall that that property is preserved by Lipschitz equivalences. The result then follows from Corollary \ref{corolnonemptyint}, Theorem \ref{main1} and an application of Pe\l czy\'nski's decomposition method.  $\Box$\\

\textbf{Remark.} Note that the compactness condition in Theorem \ref{main2} is necessary, even if have uniformity on the embeddings into $\R^n$. For example, $\Z\times\R$ is locally isometric to line segments, but $\F(\Z\times\R)$ is not isomorphic to a subspace of $\F(\R) \cong L^1$, by Naor and Schechtman's result mentioned in the introduction.


\section{$\F(X)\simeq \F(X)^2$ with low distortion}
\label{secfx2}

Let $X$ be a Banach space. By Theorem \ref{main1}, $\F(X)\simeq \F(X)^2$. 
In this Section we will show that we have the uniform bound $d_{BM}(\F(X),\F(X)\oplus_1 \F(X))\leq 4$; we will do this via an elementary construction based on metric properties of $X$.  

We start by recalling some definitions and results on quotient metric spaces which are of a more general kind than the ones presented in Section \ref{secext}. For details and more on that subject, we refer to Weaver's book \cite{W}. Let $(M,d)$ be a complete metric space, and let $\sim$ be an equivalence relation on $M$. The element of  $M/\sim$ containing $x\in M$ will be denoted by either $\tilde{x}$ or $[x]_\sim$. Define a pseudometric $\tilde{d}$ on $M/\sim$ by 
\begin{eqnarray}
\tilde{d}(\tilde{x},\tilde{y}):= \inf\left\{\sum_{j=1}^n d(x_j,y_j):n\in\N, x\sim x_1,y_j\sim x_{j+1} (j=1,...,n-1),y_n\sim y \right\}. 
\label{deftilded}
\end{eqnarray}
This pseudometric can be roughly interpreted in the following way: it is the length of the shortest discrete path from $x$ to $y$ when we are allowed to teleport between equivalent elements. An equivalent way to define $\tilde{d}$, that will be useful for further constructions, is the following: 
\begin{eqnarray}
\tilde{d}(\tilde{x},\tilde{y})= \sup\left\{|f(x)-f(y)|: f:M\rightarrow \R \mbox{ is constant in each } \tilde{z}\in \tilde{M},\|f\|_{Lip}\leq 1\right\}, 
\label{deftildedlip}
\end{eqnarray}
where $\|f\|_{Lip}$ is the Lipschitz constant of $f$. 

On $M$ we define yet the equivalence relation $\approx$ which identifies all $x,y\in M$ satisfying $\tilde{d}(\tilde{x},\tilde{y})=0$, and on $M/\approx$ we define the metric $\tilde{\tilde{d}}(\tilde{\tilde{x}},\tilde{\tilde{y}})=\tilde{d}(\tilde{x},\tilde{y})$. We define $M_\sim$, the \emph{metric quotient (or just quotient) of $M$ with respect to $\sim$,} as the completion of $M/\approx$. 
Note that, for a given complete metric space $(M,d,0)$ and an equivalence relation $\sim$ on $M$, by (\ref{deftildedlip}) there is a canonical isometric isomorphism between $Lip_{\tilde{\tilde{0}}}(M_\sim)$ and the closed subspace of $Lip_0(M)$ consisting of all functions that are constant in each class $\tilde{x}\in M/\sim$.

We recall some definitions concerning path  metric spaces.  Let $(M,d)$ be a pseudometric space, and let $\varphi: I\rightarrow M$ be a curve (that is, $I$ is an interval and $\varphi$ is continuous). The \emph{length} of $\varphi$ is $\ell(\varphi):=\sup\{\sum_{j=1}^n d(\varphi(x_{j-1}),\varphi(x_{j}))\}$, where the supremum is taken over $n\in\N$ and $x_j\in I,\, x_0<\dots<x_n$. $(M,d)$ is said to be a \emph{path metric space} if $d$ is a metric and $d(x,y)= \inf\{\ell(\varphi):\varphi$ is a curve in $M$ having endpoints $x$ and $y\}$. A \emph{minimizing geodesic} in a path metric space is any curve $\varphi:I\rightarrow M$ such that $d(\varphi(t),\varphi(s))=|t-s|$ for all $t,s\in I$; $(M,d)$ is said to be \emph{geodesic} if any two elements of $M$ are joined by a minimizing geodesic.

\begin{proposition}

Let $(M,d)$ be a path metric space. Then each metric quotient of $M$ is a path metric space. 

\label{lemmapath}
\end{proposition}

\textbf{Proof.} Fix an equivalence relation $\sim$ on $M$. Let $x,y\in M$, and for each $k\in\N$ consider pairs $(x_1^k,y_1^k),\dots,(x_{n_k}^k,y_{n_k}^k)$ of elements of $M$ such that
$$
x\sim x_1^k,y_j^k\sim x_{j+1}^k(j=1,\dots,n_k-1),y_{n_k}^k\sim y
$$ 
and $\sum_{j=1}^{n_k}d(x_j^k,y_j^k) \stackrel{k}\rightarrow \tilde{d}(\tilde{x},\tilde{y})$. Since $(M,d)$ is a path metric space, there exist, for each $k\in\N$ and $j=1,\dots,n_k$, curves $\varphi_j^k$ with endpoints $x_j^k$ and $y_j^k$, respectively, and such that 
$$
\sum_{j=1}^{n_k} \ell(\varphi_j^k) < \tilde{d}(\tilde{x},\tilde{y})+\frac1k. 
$$
Concatenating we get a curve $\tilde{\varphi}^k$ in $M/\sim$ with endpoints $\tilde{x}$ and $\tilde{y}$ satisfying $\ell(\tilde{\varphi}^k)< \tilde{d}(\tilde{x},\tilde{y})+\frac1k$. Since for any curve $\tilde{\varphi}$ in $M/\sim$ with endpoints $\tilde{x}$ and $\tilde{y}$ we have $\tilde{d}(\tilde{x},\tilde{y})\leq \ell(\tilde{\varphi})$, it follows that 
$$
\tilde{d}(\tilde{x},\tilde{y})= \inf\{\ell(\tilde{\varphi}):\tilde{\varphi}\mbox{ is a curve in $M/\sim$ having endpoints $\tilde{x}$ and $\tilde{y}$}\},  
$$ 
and then clearly the same holds for $M/\approx$ and thus for $M_\sim$.  $\Box$\\

\textbf{Remark:} In Proposition \ref{lemmapath} we cannot substitute \emph{path metric space} by \emph{geodesic metric space}: 

\begin{proposition}

There is a geodesic metric space $M$ which admits a metric quotient that is a non-geodesic path metric space. 

\end{proposition}

\textbf{Proof.} Let $e_j$ be the standard unit vectors of $\ell_1$ and consider the metric subspace of $\ell_1$ defined by $M:=\cup_{j=1}^\infty [0,1]e_j$.
Let $F:=\cup_{j=1}^\infty \{e_j\}$ and suppose that $\sim$ is the equivalence relation which collapses $F$ to a point. Note that, in this case, $M/\sim=M_\sim$, that the $M_\sim$-distance between $\tilde{0}$ and $\tilde{e_1}=F$ is $1$ and that there are minimizing geodesics with endpoints $\tilde{0}$ and $\tilde{e_1}$ going through each segment $[0,1]e_j$. 

Let $F_j:= [\frac14 + \frac{1}{2^{2+j}}, \frac34 - \frac{1}{2^{2+j}}]e_j$ be interpreted as a subset of $M$, and consider on $M$ the equivalence relation $\equiv$  that collapses each $F_j$ to a point, and the respective quotient metric space $(M/\equiv,d)$ (again, in this case we have $M_\equiv = M/
\equiv$). Then $d([\tilde{0}]_\equiv,[\tilde{e_1}]_\equiv)=\frac12$ and there are curves $\varphi_j$ in $M/\equiv$ with endpoints   $[\tilde{0}]_\equiv$ and $[\tilde{e_1}]_\equiv$ with $\ell(\varphi_j)\stackrel{j}{\rightarrow} \frac12$, even though there is no minimizing geodesic in $M/\equiv$ with endpoints    $[\tilde{0}]_\equiv$ and $[\tilde{e_1}]_\equiv$. One  verifies the path metric property for all other pairs of elements of $M/\equiv$.  $\Box$\\

\begin{lemma}

Consider $M$ a path metric space, $N$ a metric space, $f:M\rightarrow N$ and $C>0$. Then $f$ is $C$-Lipschitz if and only if it is locally $C$-Lipschitz. 

\label{lemmageod}
\end{lemma}

\textbf{Proof.} To prove the nontrivial implication, fix $\delta>0$, let $x,y\in  M$ and let $\varphi:I\rightarrow M$ be a curve with endpoints  $x$ and $y$ satisfying
$$
\ell(\varphi)<d_M(x,y)+\delta.
$$
 For each $t\in I$, there exists by hypothesis $\epsilon_t>0$ such that $f|_{\varphi(] t-\epsilon_t,t+\epsilon_t [)}$ is $C$-Lipschitz. Since $I$ is compact, there are $t_1<\dots<t_n$ such that $\cup_{j=1}^\infty  ] t_j-\epsilon_{t_j},t_j+\epsilon_{t_j} [ \supset I$. We can then easily find $\varphi$-consecutive points $z_1,\dots,z_m$ in $\varphi(I)$ satisfying
$$\begin{array}{rl}
d_N(f(x),f(y)) \!\!\! & \leq d_N(f(x),f(z_1))+d_N(f(z_1),f(z_2))+\dots +d_N(f(z_m),f(y))\\
& \leq C(d_M(x,z_1)+d_M(z_1,z_2)+\dots +d_M(z_m,y))\\
&\leq C(d_M(x,y)+\delta).
\end{array}$$
Since $\delta$ was arbitrary, the conclusion follows.  $\Box$\\

Let $(X,\|\cdot\|)$ be a Banach space. We now proceed with the construction of a pair of metric quotients of $X$, namely $X_L$ and $X_R$, which have properties useful for studying products of $Lip_0(X)$ (see Proposition \ref{FX2}). Let $\alpha,\beta: [0,+\infty)\rightarrow [0,+\infty)$ be the continuous functions defined in each $[2^m,2^{m+1}],m\in\Z$ by 
$$
\alpha(t):= \left\{ 
\begin{array}{ll}
t-2^{m-1},  &\mbox{if $2^m\leq t \leq 2^{m-1}+2^m$}; \\
2^m, &\mbox{if $2^{m-1}+2^m \leq t \leq 2^{m+1}$} 
\end{array} \right. 
$$
and
$$
\beta(t):= \left\{ 
\begin{array}{ll}
2^{m-1},  &\mbox{if $2^m\leq t \leq 2^{m-1}+2^m$};\\
t-2^m, &\mbox{if $2^{m-1}+2^m \leq t \leq 2^{m+1}$.} 
\end{array} \right.
$$
Consider the equivalence relations $\sim_L$ and $\sim_R$ on $X$ defined by
$$
x\sim_L y \Leftrightarrow x=y \mbox{ or } (x=\lambda y\mbox{ with }\lambda>0,\mbox{ and $\alpha$ is constant in }[\|x\|,\|y\|])
$$
and
$$
x\sim_R y \Leftrightarrow x=y \mbox{ or } (x=\lambda y\mbox{ with }\lambda>0,\mbox{ and $\beta$ is constant in }[\|x\|,\|y\|])
$$
and denote by $X_L = (X_L,d_L)$ and $X_R = (X_R,d_R)$ the corresponding quotient metric spaces. 

\begin{figure}
\includegraphics[height=10cm]{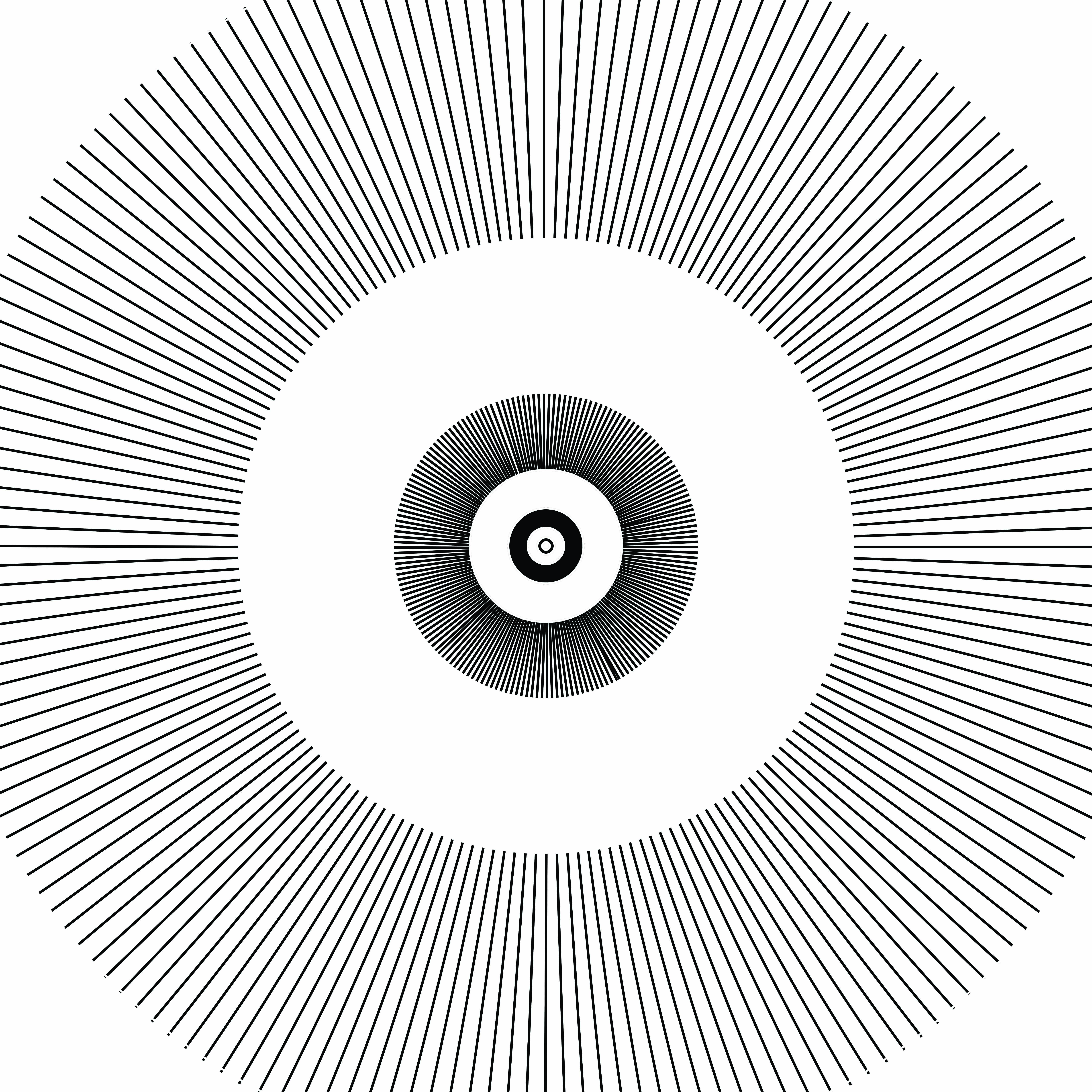}
\caption{This is how $X_L$ looks like. The represented radial segments are collapsed to points. $X_R$ looks the same, up to a factor two re-scaling.}
\end{figure}

To prove the next lemma we use Hopf-Rinow's Theorem which states that in a complete and locally compact path metric space, each pair of points are joined by a minimizing geodesic (see e.g. \cite{grombook}). 

\begin{lemma}

$X_L$ and $X_R$ are geodesic. 

\label{xrmingeod}
\end{lemma}

\textbf{Proof.} For any $x,y\in X$, note that the metric space $(span_X\{x,y\}/\sim_R,d_R)$ satisfies the conditions in Hopf-Rinow's Theorem, thus there is a minimizing geodesic $\gamma$ in $(span_X\{x,y\}/\sim_R,d_R)$ (thus also in $X_R$) with endpoints $\tilde{x}$ and $\tilde{y}$. The same argument holds for $X_L$.  $\Box$\\

\begin{lemma}

There exist onto bi-Lipschitz mappings $L:X\rightarrow X_L$ and $R:X\rightarrow X_R$ with $\|L\|_{Lip}\leq 1$, $\|L^{-1}\|_{Lip}\leq \frac43$, $\|R\|_{Lip}\leq \frac32$ and $\|R^{-1}\|_{Lip}\leq 1$.

\label{lemmaa}
\end{lemma}


\textbf{Proof.} Denote $C_m:=B_{2^{m+1}}\setminus B_{2^{m}},N\in\N$, and for each $x\in X\setminus\{0\}$ let $m_x\in \Z$ be such that $\|x\| \in C_{m_x}$. Define the bicontinuous mapping $R: X\setminus\{0\}\rightarrow X_R\setminus \{\tilde{0}\}$ by
$$
R(x):=\left(\left(\frac{1}{2} +\frac{2^{m_x}}{\|x\|}\right)x\right)^{\sim_R}.
$$
What $R$ does is to squeeze each crown $C_m$ onto the thinner crown $R(C_m)=(B_{2^{m+1}}\setminus B_{2^{m}+2^{m-1}})^{\sim_R}$. 
For $x\in X\setminus \{0\}$, let $V_x$ be a neighborhood of $x$ such that, for each $y\in V_x$, $\|x-y\|\leq 2^{m_x-1}$ and $d_R(R(x),R(y))\leq 2^{m_x-2}$. This will imply that, for any $y\in V_x$, we have $|m_x-m_y|\leq 1$, the line segment with endpoints $x$ and $y$ intercepts at most two crowns $C_m$ and a minimizing geodesic with endpoints $R(x)$ and $R(y)$ intercepts at most two crowns $R(C_m)$. We shall show that 
$$
\|x-y\|\leq d_R(R(x),R(y))\leq \frac32\|x-y\|,\,
x\in X\setminus \{0\},y\in V_x. 
$$
The fact that $X$ and $X_R$ have minimizing geodesics connecting any pair of points will allow us then to assert, by Lemma \ref{lemmageod}, that the above inequality holds without the restriction $y\in V_x$, and thus that $R$ is bi-Lipschitz, $\|R\|_{Lip}\leq \frac32$ and $\|R^{-1}\|_{Lip}\leq 1$. 

In effect, let $x,y\in X\setminus B_m$ with $\|x-y\|\leq 2^{m-1}$ and $d_R(R(x),R(y))\leq 2^{m-2}$. Assume without loss of generality that $\|x\|\leq \|y\|$. Then one of the following conditions is true:
\begin{enumerate}
\item $m_x = m_y$, and $d_R(R(x),R(y))=\|R(x)-R(y)\|$;
\item $m_x < m_y$;
\item $m_x = m_y$, and there is a minimizing geodesic with endpoints $R(x)$ and $R(y)$ passing through $R(C_{m_x-1})$. 
\end{enumerate} 
If (1) is true, then  $\frac{1}{2} +\frac{2^{m_x}}{\|y\|}\leq \frac{1}{2} +\frac{2^{m_x}}{\|x\|}$, and 
\begin{eqnarray}
\left(\frac{1}{2} +\frac{2^{m_x}}{\|y\|}\right)\|x-y\|\leq \|R(x)-R(y)\| \leq \left(\frac{1}{2} +\frac{2^{m_x}}{\|x\|}\right)\|x-y\|, 
\label{azert}
\end{eqnarray}
thus 
\begin{eqnarray}
\|x-y\|\leq d_R(R(x),R(y))\leq \frac32\|x-y\|.
\label{pqow}
\end{eqnarray} 

If (2) is true, suppose that $\|x\|<2^{m_x+1}$ (if $\|x\|=2^{m_x+1}$, then $x$ and $y$ satisfy (1))  and let $z$ be the intersection point between the line segment $[x,y]$ and $S_{2^{m_x+1}}$. Then the pairs $x,z$ and $z,y$ satisfy (1), and by (\ref{pqow}) we have 
$$
d_R(R(x),R(y))\leq d_R(R(x),R(z))+d_R(R(z),R(y))\leq \frac32 (\|x-z\|+\|z-y\|) = \frac32 \|x-y\|.
$$
Similarly, let $\tilde{z}$ be the intersection of $R(S_{2^{m_x+1}})$ with a minimizing geodesic with endpoints $R(x)$ and $R(y)$. Then $d_R(R(x),\tilde{z})=\|R(x)-\tilde{z}\|$ and $d_R(\tilde{z},R(y))=\|\tilde{z}-R(y)\|$, and thus by (\ref{pqow}) we have that
$$
\|x-y\|\leq \|x-R^{-1}(\tilde{z})\|+\|R^{-1}(\tilde{z})-y\|\leq  
d_R(R(x),\tilde{z})+d_R(\tilde{z},R(y)) =  d_R(R(x),R(y)). 
$$

For the remainder case (3), we can obtain the desired inequalities by taking a convenient point in a minimizing geodesic with endpoints $R(x)$ and $R(y)$ and reducing the problem to the case (2). \\

The Lipschitz equivalence between $X$ and $X_L$ is given by the mapping $L: X\setminus \{0\} \rightarrow X_L\setminus \{\tilde{0}\}$ defined by 
$$
L(x):=\left(\left(\frac{1}{2} +\frac{2^{m_x-1}}{\|x\|}\right)x\right)^{\sim_L},
$$
which squeezes each crown $C_m$ onto the thinner crown $L(C_m)=\left(B_{2^m+2^{m-1}}\setminus B_{2^m}\right)^{\sim_L}$. 
To show this and obtain the Lipschitz constants, one simply must follow the same steps taken to do that for $R$. The only difference will appear when getting to (\ref{azert}), which will read
$$
\left(\frac{1}{2} +\frac{2^{m_x-1}}{\|y\|}\right)\|x-y\|\leq \|L(x)-L(y)\| \leq \left(\frac{1}{2} +\frac{2^{m_x-1}}{\|x\|}\right)\|x-y\|,
$$
and thus (\ref{pqow}) will read 
$$
\frac34 \|x-y\|\leq d_L(L(x),L(y))\leq \|x-y\|.
$$
Following analogous steps, we get to the conclusion.  $\Box$\\

\begin{proposition}

Let $X$ be a Banach space. Then $d_{BM}(\F(X),\F(X)\oplus_1 \F(X))\leq 4$. 

\label{FX2}
\end{proposition}

\textbf{Proof.} Recall that, by (\ref{deftildedlip}), $Lip_0(X_L)\cong Y_L$ and $Lip_0(X_R)\cong Y_R$, where $Y_L$ and $Y_R$ are the closed subspaces of $Lip_0(X)$ defined by
$$
Y_L:= \{f\in Lip_0(X): \mbox{$f$ is constant in each equivalence class of $X_L$}\}
$$
and
$$
Y_R:= \{f\in Lip_0(X): \mbox{$f$ is constant in each equivalence class of $X_R$}\}.
$$
Let $\Phi:Y_L \oplus_\infty Y_R \rightarrow Lip_0(X)$ be defined by  $\Phi(f,g) :=f+g$. Then $\Phi$ is linear with norm $\|\Phi\|\leq2$. Moreover, $\Phi$ admits an inverse defined by
$$
(\Phi^{-1}h)(x)=\left(\frac{\alpha ( \|h(x)\| )}{\|h(x)\| } h(x),\frac{\beta( \|h(x)\| )}{\|h(x)\| } h(x)\right).
$$
Since the functions $x\mapsto \frac{\alpha ( \|x\| )}{\|x\| } x$ and $x\mapsto \frac{\beta ( \|x\| )}{\|x\| } x$ are 1-Lipschitz, it follows that $\|\Phi^{-1}\|\leq 1$. Now from Lemma \ref{lemmaa} it follows that  there is an isomorphism $\Psi$ from $Lip_0(X)\oplus_\infty Lip_0(X)$ onto $Y_L\oplus_\infty Y_R$ satisfying $\|\Psi\|.\|\Psi^{-1}\|\leq \frac43.\frac32=2$. Then $\Phi\circ\Psi$ is an isomorphism from $Lip_0(X)\oplus_\infty Lip_0(X)$ onto $Lip_0(X)$ satisfying $\|\Phi\circ\Psi\|.\|(\Phi\circ\Psi)^{-1}\|\leq 4$. Since $\Phi$ and $\Psi$ are pointwise-to-pointwise continuous, $\Phi\circ\Psi$ induces an isomorphism $T: \F(X)\rightarrow \F(X)\oplus_1\F(X)$ satisfying $T^*=\Phi\circ\Psi$ and $\|T\|.\|T^{-1}\|\leq 4$.  $\Box$\\


\section{On Lipschitz-free spaces over unions of metric spaces}
\label{secunion}

In this section we will provide a couple of formulas for computing the Lipschitz-free space over certain unions of metric spaces from the Lipschitz-free spaces of the original metric spaces, provided that we have an \emph{orthogonal} placement of the metric spaces involved (in a sense to be made precise). The first one generalizes in particular a first step taken in this direction by Godard \cite{G} (see Proposition \ref{propgod} below). The idea of studying the behavior of the Lipschitz-free spaces with respect to unions is motivated, for example, by the problem of characterizing the metric spaces such that the corresponding Lipschitz-free spaces admit isomorphic embeddings into $L^1$. We start by studying how taking certain \emph{orthogonal} unions of metric spaces have effect on the corresponding Lipschitz-free spaces:

\begin{proposition}

Suppose that $M=\cup_{\gamma\in\Gamma} M_\gamma$ is a metric space with metric $d$, and suppose that there exists $0\in M$ satisfying
\begin{enumerate}
\item $M_\gamma\cap M_\eta = \{0\}$ for $\gamma\neq \eta$, and
\item \emph{(orthogonality)} there exists $C\geq 1$ such that, for all $\gamma\neq \eta,\, x\in M_\gamma$ and $y\in M_\eta$, $d(x,0)+d(y,0)\leq C\,d(x,y)$. 
\end{enumerate}
Then
\begin{eqnarray*}
\F \left(\bigcup_{\gamma\in \Gamma}M_\gamma\right) \simeq \left(\sum_{\gamma\in \Gamma}\F (M_\gamma)\right)_{\ell_1}
\label{isounion}
\end{eqnarray*}
with distortion bounded by $C$.
 
\label{propbasic2}
\end{proposition}

\textbf{Proof.} We can assume $0 = 0_M = 0_{M_\gamma}$, for all $\gamma$. Consider 
\begin{eqnarray}
\Phi: (f_\gamma)\in \left(\sum_{\gamma\in\Gamma}Lip_0(M_\gamma)\right)_{\ell_\infty}\mapsto h\in Lip_0(\cup_{\gamma\in\Gamma}M_\gamma), 
\label{concat}
\end{eqnarray}
where $h|_{M_\gamma}=f_\gamma$. $h=\Phi (f_\gamma)$ is well defined since, for each $x\in M_\zeta$ and $y\in M_\eta$ with $\zeta\neq\eta$,
\begin{align*}
|\Phi((f_\gamma))(x)-\Phi((f_\gamma))(y)| &=|f_\zeta(x)-f_\eta(y)|\leq \|f_\zeta\|d(x,0) +\|f_{\eta}\|d(y,0)\\
& \leq C\,max\{\|f_\zeta\|,\|f_\eta\|\}d(x,y)
\leq C\|(f_\gamma)\|d(x,y). 
\end{align*}
Since $\Phi$ is linear, we have from the inequality above that 
\begin{eqnarray}
\|\Phi\|\leq max\{1,C\}.
\label{max1C}
\end{eqnarray} 
It is clear that $\Phi$ is surjective, and that $\|\Phi^{-1}\| \leq 1$. $\Phi$ is weak$^*$-to-weak$^*$ continuous, thus it is the adjoint of an isomorphism from $\F(\cup_{\gamma\in\Gamma}M_\gamma)$ onto $\left(\sum_{\gamma\in\Gamma}\F(M_\gamma)\right)_{\ell_1}$. It is clear from (\ref{max1C}) that the distortion of $\Phi$ is bounded by $C$.  $\Box$\\

Note that, if we exclude the orthogonality condition, the concatenating operator $\Phi$ defined in (\ref{concat}) does not need to take values in $Lip_0(M)$. We do not know if it is possible to remove this condition, even for finite unions. In fact, in particular we do not now the answer to the following apparently simple question: \\

\textbf{Problem.} Consider the set $Cusp:=\{(x,0):x\geq 0\}\cup\{(x,x^2):x\geq 0\}$ endowed with the euclidean metric. Is $\F(Cusp)$ isomorphic to a subspace of $L^1$?\\

Note that $Cusp$ is not Lipschitz equivalent to the real line, so we cannot appeal to Proposition \ref{propgk} to have a positive answer. A negative answer would imply that there is no isomorphic copy of $\F(\R^2)$ in $L^1$, which as we mentioned is one of the main results from \cite{NS}.\\

The mentioned result by Godard can be seen as a particular case of Proposition \ref{propbasic2}:

\begin{proposition}[Godard \cite{G}, Proposition 5.1]

Let $\Gamma$ be a set with a distiguished point $0$, denote $\Gamma^\ast:=\Gamma\setminus \{0\}$, and let $M=\cup_{\gamma\in\Gamma}M_\gamma$ be a metric space with metric $d$. Suppose that there exist $A,B>0$ such that $A\leq d(x,y) \leq B$ whenever $x$ and $y$ belong to different $M_\gamma$'s. Then 
$$
\F(M) \simeq \left(\sum_{\gamma\in\Gamma} \F(M_\gamma)\right)_{\ell_1} \oplus \ell_1(\Gamma^\ast).
$$ 

\label{propgod}
\end{proposition}

\textbf{(Alternative) proof.}  Assume that $A\leq 1\leq B$. Fix a point $p\in M_0$, and for each $\gamma\neq 0$ fix $0_\gamma\in M_\gamma$ and define $\Phi_\gamma = (\Phi_\gamma^1,\Phi_\gamma^2): Lip_p(M_\gamma\cup\{p\})\rightarrow Lip_{0_\gamma}(M_\gamma)\oplus_\infty \R$ by 
$$(\Phi_\gamma^1(f),\Phi_\gamma^2(f)):=(f-f(0_\gamma),f(0_\gamma)).$$
It is easily seen that $\|\Phi_\gamma\|\leq B$. $\Phi_\gamma$ admits an inverse $\Phi_\gamma^{-1}$ which is defined by $\Phi_\gamma^{-1}(f,r)(x)=f(x)+r$, if $x\in M_\gamma$ and $\Phi_\gamma^{-1}(f,r)(p)=0$. We have the bound $\|\Phi_\gamma^{-1}\|\leq \frac{B+1}{A}$; to see this, let $f\in B_{Lip_{0_\gamma}(M_\gamma)}$ and $r\in\R$ with $|r|\leq 1$, and let $x\in M_\gamma$. Then 
$$
\frac{| \Phi_\gamma^{-1}(f,r)(x)- \Phi_\gamma^{-1}(f,r)(p) |}{d(x,p)} \leq \frac{| f(x) + r |}{d(x,p)}\leq 
\frac{|f(x)| + |r|}{A}\leq \frac{B+1}{A}.  
$$
It follows that $\Phi:Lip_p(M_0)\oplus_\infty \left(\sum_{\gamma\in\Gamma^*} Lip_p(M_\gamma \cup\{p\})\right)_{\ell_\infty}\rightarrow \left(\sum_{\gamma\in\Gamma}Lip_{0_\gamma}(M_\gamma)\right)_{\ell_\infty} \oplus_\infty \ell_\infty(\Gamma^*)$ defined by 
$$
\Phi(f,(f_\gamma)_{\gamma\in\Gamma^*}):=((f,(\Phi_\gamma^1 (f_\gamma))_{\gamma\in\Gamma^*}),(\Phi_\gamma^2 (f_\gamma))_{\gamma\in\Gamma^*} )
$$
is a pointwise-to-pointwise continuous isomorphism with $\|\Phi\|.\|\Phi^{-1}\|\leq \frac{B(B+1)}{A}$. It is therefore the adjoint of an isomorphism 
\begin{eqnarray}\F(M_0)\oplus_1 \left(\sum_{\gamma\in\Gamma^*} \F(M_\gamma \cup\{p\})\right)_{\ell_1}
\simeq\left(\sum_{\gamma\in\Gamma}\F(M_\gamma)\right)_{\ell_1} \oplus_1 \ell_1(\Gamma^*).
\label{dfv}
\end{eqnarray}
Now the spaces $M_0,(M_\gamma\cup\{p\}),\gamma\in\Gamma^*$ satisfy the orthogonality condition (2) of Proposition \ref{propbasic2}, thus $\F(M) = \F(M_0\cup (\cup_{\gamma\in\Gamma^*}   (M_\gamma\cup\{p\}))  \simeq \F(M_0)\oplus_1 \left(\sum_{\gamma\in\Gamma^*} \F(M_\gamma \cup\{p\})\right)_{\ell_1}$, and the conclusion follows by (\ref{dfv}).  $\Box$\\

Note that, using Lemma \ref{lemmaextFM}, we can also have some control when the intersection of the metric spaces involved is nontrivial. We present a version for the union of two metric spaces.

\begin{proposition}

Let $M\cup N$ be a metric space with metric $d$, suppose that $F:= M\cap N$ is closed and nonempty, and choose a base point $0$ in $F$. Assume that
\begin{enumerate}
\item there exists $E$ in $Ext_0^{pt}(F,M\cup N)$ satisfying $\|E\|\leq C$, and
\item \emph{(orthogonality)} there exists $C\geq 1$ such that, for all $x\in M$ and $y\in N$, $d(x,F)+d(y,F)\leq C\,d(x,y)$. 
\end{enumerate}
Then
\begin{eqnarray*}
\F (M\cup N) \simeq \F (M/F) \oplus_1 \F (N/F) \oplus_1 \F (F)
\label{formunion}
\end{eqnarray*}
with distortion bounded by $C(\|E\|+1)^2$. 

\label{UNION2}
\end{proposition}


\textbf{Proof.} Denote by $\tilde{d}$ the metric of  $(M\cup N)/F$, and note that $M/F$ and $N/F$ are subsets of  $(M\cup N)/F$. For each $\tilde x\in M/F$ and each $\tilde y\in N/F$, 
\begin{align*}
\tilde d(\tilde x,\tilde 0)+\tilde d(\tilde y,\tilde 0) &= min\{d(x,0),d(x,F)\}+min\{d(y,0),d(y,F)\}\\
&=d(x,F)+d(y,F) \leq C\, min\{d(x,y),d(x,F)+d(y,F)\} = C\,\tilde d(\tilde x,\tilde y),
\end{align*}
and thus we can apply Proposition \ref{propbasic2} and get that  $d_{BM}(\F(M\cup N/F),\F(M/F)\oplus_1 \F(N/F))\leq C$.
Since $\F(M\cup N)$ is $(\|E\|+1)^2$-isomorphic to $\F((M\cup N)/F)\oplus_1 \F(F)$ by Lemma \ref{lemmaextFM}, the result follows.  $\Box$\\

\textbf{Acknowledgements} This research was supported by CAPES, grant BEX 10057/12-9.

I am grateful for having been warmly welcomed by the \'Equipe d'Analyse Fonctionnelle during my postdoctoral period at the IMJ. 
I am especially grateful to Gilles Godefroy for many fruitful discussions which lead to several improvements of this work. I thank Florent Baudier for driving my attention to \cite{LanPla}. I thank Bruno Braga, Marek C\'uth and Michal Doucha for reading my original manuscript and pointing out flaws to be corrected.  Illustration by Gabriela Kaufmann Sacchetto.

\end{document}